\documentclass[a4paper,10pt]{article}

\usepackage{amsmath, amsfonts, amssymb, latexsym, mathrsfs, graphicx}

\usepackage[all]{xy}
\newdir{ >}{{}*!/-5pt/@{>}}
\newdir{  >}{{}*!/-7pt/@{>}}

\newtheorem{remark}{Remark}

\newcommand{\mmap}[3]{\ar@/#2/[#1]|*=0{\rotatebox{#3}{$\scriptsize |$}}} 
\newcommand{\picar}[3]{\ar@{}[#1]|*{\rotatebox{#2}{$\scriptsize #3$}}} 

\newcommand{\ba}{{\mathscr{BA}}}
\newcommand{\sa}{{\mathscr{A}}}

\newcommand{\scrb}{{\mathscr{B}}}
\newcommand{\sv}{{\mathscr{V}}}
\newcommand{\ssc}{{\mathscr{C}}}
\newcommand{\oop}{\operatorname{op}}
\newcommand{\op}{^{\oop}}
\newcommand{\set}{\mathbf{Set}}
\newcommand{\vect}{\mathbf{Vect}_{k}}

\begin{document}

\begin{center}
{\underline{{\Large Note On Analytic Functors As Fourier Transforms. }}}
\newline

{{\small Brian J. DAY}}

{{\small April 19, 2009.}}
\end{center}

\noindent \emph{\small \underline{Abstract:} Several notions of ``analytic'' functor introduced recently in the literature fit into the graphic fourier transform context presented in [D].}
\newline

\begin{center}
\begin{tabular}{c}\hline
\phantom{a}\phantom{a}\phantom{a}\phantom{a}\phantom{a}\phantom{a}\phantom{a}\phantom{a}\phantom{a}\phantom{a}\phantom{a}\phantom{a}\phantom{a}\phantom{a}\phantom{a}\phantom{a}
\end{tabular}
\end{center}

Various concepts of ``analytic'' functor are well characterized in different places in the literature. Here we want to mention explicitly how two ideas, introduced in [AV] and [FGHW], can also be viewed from [D].
However, we don't offer further characterization results here.
\newline

\noindent\underline{Example [FGHW]:}
\newline

Let $N: \ba \to [\sa\op,\set]$ be the canonical ``inclusion'' functor from the (monoidal) groupoid $\ba$ constructed in [FGHW] from the small $\set$-category $\sa$; i.e. $\ba$ is the groupoid of all isomorphisms in the free finite-coproduct completion of $\sa$ in $[\sa\op,\set]$.
Then, assuming (here) that $\ssc$ is small, the functor

$$\exists_{N} : [\ba,[\ssc\op,\set]] \to [[\sa\op,\set],[\ssc\op,\set]]$$

\noindent which is precisely the process of left Kan extension along N, is conservative (because $\ba$ is a groupoid) and tensor product preserving (because N preserves finite coproducts).
This is then consistent with [D] for $\sv = [\ssc\op,\set]$ (cartesian monoidal).
\newline

\noindent\underline{Example [AV]:}
\newline

In [AV] Remark 4.5, the authors complain (justly) that many of their ``analytic'' functors for $\sv = \vect$ are not $k$-linear.
But this is not too serious a matter because the set-up in [D] permits a $k$-linearization (this is not a tautology, but merely an adjunction).
Thus the ``ordinary'' kernel discussed in [AV] \S 4, namely

$$K : \scrb \times \sv_{0} \longrightarrow \sv_{0}, \,\, K(n,X) = \otimes^{n}X ,$$

\noindent yields the corresponding (multiplicative) $\sv$-kernel 

$$\xymatrix{E : k_{\ast}\scrb\otimes k_{\ast}\sv_{0} \ar[rr]^-{\cong}&& k_{\ast}(\scrb\times\sv_{0}) \ar[rr]^-{k_{\ast}K} && k_{\ast}\sv_{0} \ar[rr]^-{can.} && \sv }$$

\noindent where $k_{\ast}$ denotes the free-$\vect$-structure functor, and $\sv_{0}$ is the ordinary category underlying $\sv$ (here $k_{\ast}\sv_{0}$ has the comonoidal structure directly induced by that on $\sv_{0}$).
Then the $\sv$-functor  

$$\overline{E} : [k_{\ast}\scrb,\sv] \longrightarrow [k_{\ast}\sv_{0},\sv] ,$$ 

\noindent is conservative (since $\scrb$ is a groupoid) and tensor product preserving (since $E$ is multiplicative).
The Fourier transforms $\overline{E}(f)$ can thus be viewed as either $k$-linear ``$E$-analytic'' functors

$$k_{\ast}\sv_{0} \longrightarrow \sv ,$$

\noindent or just ordinary ``[AV]-analytic'' functors 

$$\sv_{0} \longrightarrow \sv_{0},$$
in the sense of [AV] Definition 4.1. 
Then the considerations of [D] Section 1.3 apply.

\begin{remark}
The term ``analytic'' functor seems quite appropriate in such cases.
\end{remark}

\noindent\underline{Example [D]:}
\newline

A type of ``quantum category'' example evolves from any $\sv$-promonoidal category $(\sa,p,j)$.
Namely, the left ``Cayley'' functor 

$$\overline{K} : [\sa,\sv] \longrightarrow [\sa\op\otimes\sa,\sv]$$

\noindent given by

$$\overline{K}(f)(A,B) = \int\limits^{X} p(X,A,B) \otimes f(X),$$

\noindent the $\sv$-kernel functor 

$$K : \sa\op\otimes\sa\op\otimes\sa\longrightarrow\sv$$

\noindent here being just the promultiplication $p$.
This $\overline{K}$ is both conservative and tensor preserving, where $[\sa,\sv]$ has the convolution structure and $[\sa\op\otimes\sa,\sv]$ has the tensor product defined by bimodule composition.
Thus $\overline{K}$ qualifies as a ``Fourier transformation'' [D].

\pagebreak

\begin{center}
{\underline{{\Large References.}}}
\end{center}

\begin{center}
\begin{tabular}{l l}
\textrm{[AV]} & J. Adamek and J. Velebil, ``Analytic functors and weak pullbacks'', \\
 & Theory Appl. Categories, 21(11), (2008) 191-209.\\
 & \\
\textrm{[D]} & B. J. Day, ``Monoidal functor categories and graphic Fourier transforms'', \\
 & arXiv:mathQA/0612496v1, 18 Dec. 2006.  \\
 & \\
\textrm{[FGHW]} & M. Fiore, N. Gambino, M. Hyland and G. Winskel, ``The cartesian closed \\ 
 & category of generalized species of structures'', London Math. Soc. (2007), 1-18. \\
 & \\
 & \\
 & \\
 & \\
\end{tabular}

\small{Mathematics Dept., Faculty of Science, Macquarie University, NSW 2109, Australia.}
\newline

Any replies are welcome through Tom Booker (thomas.booker@students.mq.edu.au), who kindly typed the manuscript.

\end{center}

\end{document}